\documentclass{article}
\usepackage{amssymb}
\usepackage{amsfonts}
\usepackage{amsmath}

\begin{document}

\begin{center}
{\Large \ Some examples of division symbol algebras of degree }${\Large 3}$%
{\Large \ and }${\Large 5}${\Large \ \\[0pt]
}

\begin{equation*}
\end{equation*}%
Cristina FLAUT and Diana SAVIN%
\begin{equation*}
\end{equation*}
\end{center}

\textbf{Abstract. }{\small In this paper we} {\small provide an algorithm to
compute the product between two elements in a symbol algebra of degree }$n$ 
{\small and\ we find an octonion non-division algebra in a symbol algebra of
degree three. Starting from this last idea, we try to find an answer to the
question if there are division  symbol  algebras of degree three. The answer
\ is positive and we provide, using MAGMA software, some examples of division symbol algebras of
degree \ }$3$ {\small and of degree} $5.$ {\small Moreover, we will give
some interesting applications of the symbol algebras in number theory. }

\bigskip 
\begin{equation*}
\end{equation*}

\bigskip


.\newline
\begin{equation*}
\end{equation*}

\textbf{0. Preliminaries}%
\begin{equation*}
\end{equation*}

Let \ $K$ be a field which contains a primitive $n$-th root of unity, with $%
n $ an arbitrary positive integer such that $char(K)$ does not divide $n.\ $%
Let $K^{\ast }=K\backslash \{0\},$ $a,b$ $\in K^{\ast }$ and let $S$ be an
algebra over $K$ generated by the elements $x$ and $y$ where%
\begin{equation}
x^{n}=a,y^{n}=b,yx=\xi xy.  \tag{1.}
\end{equation}%
where $\xi $ is a primitive root of order $n$ of unity. This algebra is
called a \textit{symbol algebra }(also known as a \textit{power norm residue
algebra}) and it is denoted by $\left( \frac{a,~b}{K,\xi }\right) .$ In [Mi;
71], J. Milnor calls this algebra "the symbol algebra" because of its
connection with the $K$-theory and with the Steinberg symbol. Symbol
algebras generalize the quaternion algebras (for $n=2$). Quaternion algebras
and symbol algebras are important not only for the theory of associative
algebras. They have many applications, some of them being studied \ by the
authors of this article: in number theory ([Sa, Fl, Ci; 09], [Mil; 10]), in
representation theory ([Fl, Sa; 13]) or in analysis and mecanics ([Ja, Ya;
13]).

In this paper, we will study the symbol algebras from two points of view:
from the theory of associative algebras and from number theory.

The study of symbol algebras of degree $n$ involves very complicated
calculations and, usually, can be hard to multiply two elements or find \
examples for some notions. In this paper, we will provide an easy algorithm
which allows us quickly computing of two elements in a symbol algebra. Since
for $n=2$ the quaternion algebras are symbol algebras, a natural question
is: for $n=3,$what is the connexion between the octonion algebras, algebras
of dimension $8,$ and symbol algebras of degree $3.$ The answer is that we
always can find an octonion non-division algebra in a symbol algebra of
degree three. Starting from this idea and from results obtained in the paper
[Fla; 12], in which, using the associated trace form for a symbol algebra,
the author studied some properties of \ such objects \ and gave some
conditions for a symbol algebra to be with division or not only for $n=4k+2$
(and not for $n\in \{3,5\}),$ we intend to find examples of division symbol
algebras of degree $3$ and degree $5$, proving that such algebras can be
with division. Since such an example is not easy to provide, we will use
MAGMA software.

In the following, we will recall some general properties and definitions.

Let $A$ be a finite dimensional unitary algebra over a field $\ K$ with a 
\textit{scalar} \textit{involution} $\,\,\,\,\overline{\phantom{x}}%
:A\rightarrow A,a\rightarrow \overline{a},\,\,$ i.e. a linear map satisfying
the following relations:$\,\,\,\,\,\overline{ab}=\overline{b}\overline{a},\,%
\overline{\overline{a}}=a,\,\,$and $a+\overline{a},a\overline{a}\in K\cdot
1\ $for all $a,b\in A.$The element $\,\overline{a}$ is called the \textit{%
conjugate} of the element $a,$ the linear form$\,\,t:A\rightarrow
K\,,\,\,t\left( a\right) =a+\overline{a}$ \ and the quadratic form $%
n:A\rightarrow K,\,\,n\left( a\right) =a\overline{a}$ are called the \textit{%
trace} and the \textit{norm \ }of \ the element $a.$

Let$\,\,\,\gamma \in K$ \thinspace be a fixed non-zero element. We define
the following algebra multiplication on the vector space 
\begin{equation*}
A\oplus A:\left( a_{1},a_{2}\right) \left( b_{1},b_{2}\right) =\left(
a_{1}b_{1}+\gamma \overline{b_{2}}a_{2},a_{2}\overline{b_{1}}%
+b_{2}a_{1}\right) .
\end{equation*}%
\newline
We obtain an algebra structure over $A\oplus A,$ denoted by $\left( A,\gamma
\right) $ and called the \textit{algebra obtained from }$A$\textit{\ by the
Cayley-Dickson process.} $\,$We have $\dim \left( A,\gamma \right) =2\dim A$.

Let $x\in \left( A,\gamma \right) $, $x=\left( a_{1},a_{2}\right) $. The map 
\begin{equation*}
\,\,\,\overline{\phantom{x}}:\left( A,\gamma \right) \rightarrow \left(
A,\gamma \right) \,,\,\,x\rightarrow \bar{x}\,=\left( \overline{a}_{1},\text{%
-}a_{2}\right) ,
\end{equation*}%
\newline
is a scalar involution of the algebra $\left( A,\gamma \right) $, extending
the involution $\overline{\phantom{x}}\,\,\,$of the algebra $A.$

\thinspace If we take $A=K$ \thinspace and apply this process $t$ times, $%
t\geq 1,\,\,$we obtain an algebra over $K,\,\,A_{t}=\left( \frac{\alpha
_{1},...,\alpha _{t}}{K}\right) .$

By induction in this algebra, the set $\{1,e_{2},...,e_{n}\},n=2^{t},$
generates a basis with the properties: $e_{i}^{2}=\gamma _{i}1,\,\,\gamma
_{i}\in K,\gamma _{i}\neq 0,\,\,i=2,...,n$ \ and\newline
$e_{i}e_{j}=-e_{j}e_{i}=\beta _{ij}e_{k},\,\,\beta _{ij}\in K,\,\,\beta
_{ij}\neq 0,i\neq j,i,j=\,\,2,...n,$ \ $\ \beta _{ij}$ and $e_{k}$ being
uniquely determined by $e_{i}$ and $e_{j}.$ For $n=2,$ we obtain the
quaternion algebra, for $n=3,$ we obtain the octonion algebra, etc.

For details about the Cayley-Dickson process, the reader is referred to $%
\left[ \text{Sc; 66}\right] $ and [Sc; 54].

If an algebra $A$ is finite-dimensional, then it is a division algebra if
and only if $A$ does not contain zero divisors. (See [Sc;66])

A central simple algebra $A$ over a field $K$ is called \textit{split by} $L$
(where $L$ is a field containing $K),$ if $A\otimes _{K}L$ is a matrix
algebra over $L.$ We also can say that its class is in the Brauer group Br$%
\left( K\right) $. (see [Ir, Ro; 92]) $L$ is called a \textit{splitting field%
} for $A.$\medskip \smallskip \newline
\textbf{Theorem 1.1.} ([Gi, Sz; 06]) \textit{Let }$K$\textit{\ be a field
such that }$\xi \in K,\,\,\xi ^{n}=1,\xi $\textit{\ is a primitive root, and
let }$\alpha ,\beta \in K^{\ast }.$\textit{\ Then the following statements
are equivalent:}\newline
\textit{i) The cyclic algebra }$A=\left( \frac{\alpha ,\beta }{K,\xi }%
\right) $\textit{\ is split.}\newline
\textit{ii) The element }$\beta $\textit{\ is a norm from the extension }$%
K\subseteq K(\sqrt[n]{\alpha }).\medskip $ \smallskip \newline
\textbf{Theorem 1.2.} ([Ir, Ro; 92]) \textit{Let }$l$\textit{\ be a natural
number, }$l\geq 3$\textit{\ and }$\xi $\textit{\ be a primitive root of the
unity of\thinspace \thinspace }$l$-\textit{order. If }$p$\textit{\ is a
prime natural number, }$l$\textit{\ is not divisible with }$p$\textit{\ and }%
$f$\textit{\ is the smallest positive integer such that }$p^{f}\equiv 1$%
\textit{\ mod }$l$\textit{, then we have} 
\begin{equation*}
p\mathbb{Z}[\xi ]=P_{1}P_{2}....P_{r},
\end{equation*}%
\textit{where }$r=\frac{%
\begin{array}{c}
\\ 
\varphi \left( l\right)%
\end{array}%
}{%
\begin{array}{c}
f%
\end{array}%
},\varphi $\textit{\ is the Euler's function and }$P_{j},\,j=1,...,r$\textit{%
\ are different prime ideals in the ring }$\mathbb{Z}[\xi ].$\smallskip 
\newline
\textbf{Theorem 1.3.} ([Lem; 00]) \textit{Let }$\xi $\textit{\ be a
primitive root of the unity of }$l-$\textit{order, where }$l$\textit{\ is a
prime natural number and let} $A$ \textit{be the ring of integers of the
Kummer field} $Q(\xi ,\sqrt[l]{\mu })$ . \textit{A prime ideal} $P$\textit{\
in the ring} $\mathbb{Z}[\xi ]$\textit{\ is in }$A$\textit{\ in one of the
situations:}

\textit{i) It is equal with the }$l-$\textit{power of a prime ideal from} $%
A, $ \textit{if the} $l-$\textit{power character }$\left( \frac{\mu }{P}%
\right) _{l}=0;$

\textit{ii) It is a prime ideal in} $A$, \textit{if} $\left( \frac{\mu }{P}%
\right) _{l}=$\textit{\ a rot of order }$l$\textit{\ of unity, different
from }$1$. \newline
\textit{iii) It decomposes in }$l$\textit{\ different prime ideals from} $A$%
, \textit{if} $\left( \frac{\mu }{P}\right) _{l}=1.$\newline
\smallskip \newline
\textbf{Theorem 1.4.} ([Mil; 10]) \textit{Let} $K$ \textit{be a finite field
and let} $Br\left( K\right) $ \textit{be} \textit{the Brauer group of } $K.$ 
\textit{Then} $Br\left( K\right) =0.$\newline
\smallskip \newline
\textbf{Theorem 1.5.} ([Al, Io; 84]) \textit{Let} $L/K$ \textit{be an
extension of finite fields. Then the norm function} $N_{L/K}:L^{\ast
}\rightarrow K^{\ast }$ \textit{is surjective }. \newline
\smallskip \newline
\textbf{Remark 1.6.} ([Led; 05]) Let $K$ be a field of characteristic $\neq
p,p$ prime, and let $\xi $$\in $$K$ be a primitive root of unity of order $p 
$\textit{. }For $a,b$$\in $$K^{\ast },$ the symbol algebra of degree $p$
denoted by $A=$ $\left( \frac{a,~b}{K,\xi }\right) $ is either split or a
division algebra\textit{. }From here, in hole this paper, we will use the
notion "no-division" instead of "split", for all symbol algebras of degree $%
n,$ with $n$ a prime number$.$

\begin{equation*}
\end{equation*}

\textbf{2. Multiplication table for symbol algebras}%
\begin{equation*}
\end{equation*}

In [Ba; 09], the author described \ how we can multiply the basis vectors in
all algebras obtained by the Cayley-Dickson process. Since the quaternion
algebra is an algebra obtained by this process and in the same time is a
particular case of symbol algebras, we use some ideas given in this paper
for multiplication of two symbol elements.

\bigskip

\textbf{Case} \thinspace $n=3.$

Let $S$ be a symbol algebra of degree three with the basis 
\begin{equation}
B=\{1,x,x^{2},y,y^{2},xy,x^{2}y^{2},x^{2}y,xy^{2}\}.  \tag{2.1.}
\end{equation}

\bigskip

\textbf{Remark 2.1.} The elements from the basis $B$ will be denoted such as
follows:\newline
$%
e_{0}=1,e_{1}=y,e_{2}=y^{2},e_{3}=x,e_{4}=xy,e_{5}=xy^{2},e_{6}=x^{2},e_{7}=x^{2}y,e_{8}=x^{2}y^{2}. 
$ If we use the lexicographic order for the monomials $x^{i}y^{j},$ we have
that $x^{i}y^{j}\geq x^{p}y^{q}$ if and only if $i\geq p$ or $i=p$ and $%
j\geq q.$ Therefore the elements from the basis $B$ are lexicographic
ordered.\medskip

\textbf{Remark 2.2. }\newline
If we write \newline
$4=1\cdot 3+1=011_{3}\rightarrow e_{4}=x^{1}y^{1}$\newline
$5=1\cdot 3+2=012_{3}\rightarrow e_{5}=x^{1}y^{2}$\newline
$6=2\cdot 3=020_{3}\rightarrow e_{6}=x^{2}$\newline
$7=2\cdot 3+1=021_{3}\rightarrow e_{7}=x^{2}y^{1}$\newline
$8=2\cdot 3+2=022_{3}\rightarrow e_{8}=x^{2}y^{2},$ where $0ij_{3}=i\cdot
3+j=k,i,j\in \{1,2\}$ is the ternary decomposition of the natural number $%
k\in \{4,5,6,7,8\},$ it results that $e_{k}=x^{i}y^{j},$with $k=i\cdot
3+j=0ij_{3}.$

If we compute two elements of the basis $B,$ we obtain%
\begin{equation}
e_{i}e_{j}=\alpha \left( i,j\right) e_{i\circ j},  \tag{2.2.}
\end{equation}%
where $\alpha \left( i,j\right) $ is a function $\alpha :\mathbb{Z}%
_{3}^{3}\times \mathbb{Z}_{3}^{3}\rightarrow K$ and $i\circ j$ represents
the "sum" of $i$ and $j$ in the group $\mathbb{Z}_{3}^{3}$ (here $i$ and $j$
are in the ternary forms!)$.$ Indeed, this last sentence results from
relation $\left( 1\right) .$ Therefore $\ e_{7}e_{8}=\alpha \left(
7,8\right) e_{3},$ since $021_{3}+022_{3}=010_{3}\rightarrow 3.\medskip $

\textbf{General case\medskip }

Using the above notations, a basis in a symbol algebra of degree $n$ is on
the form 
\begin{equation}
B=\{x^{i}y^{j}\text{ }/\text{ }0\leq i<n,0\leq j<n\}.  \tag{2.3.}
\end{equation}

The elements from the basis $B$ are lexicographic ordered, as in Remark 2.1.
We denote an element from the basis $B$ given by $\left( 2.3.\right) ,\ $\
with $e_{k}=x^{i}y^{j},0\leq k<n,$ such that $k=0ij_{n}=i\cdot n+j,$ where $%
0ij_{n}$ is the $n-$ary decomposition of the natural number $k.$ Then, using
relation $\left( 1\right) ,$ if we compute two elements of the basis $B,$ we
obtain%
\begin{equation}
e_{i}e_{j}=\alpha \left( i,j\right) e_{i\circ j},  \tag{2.4.}
\end{equation}%
where $\alpha \left( i,j\right) $ is a function $\alpha :\mathbb{Z}%
_{n}^{n}\times \mathbb{Z}_{n}^{n}\rightarrow K$ and $i\circ j$ represents
the "sum" of $i$ and $j$ in the group $\mathbb{Z}_{n}^{n}~$(with $i$ and $j$
are in the $n-$nary forms!)$.$

If $e_{i}=x^{i_{1}}y^{i_{2}}$ and $e_{j}=x^{j_{1}}y^{j_{2}},$ then we have 
\newline
\begin{equation}
e_{i}e_{j}=x^{i_{1}}y^{i_{2}}x^{j_{1}}y^{j_{2}}\rightarrow x^{i_{1}}\underset%
{i_{2}}{\underbrace{yy...y}}~\underset{j_{1}}{\underbrace{xx...x}}~y^{j_{2}}%
\overset{\omega }{\rightarrow }x^{i_{1}}\underset{i_{2}}{\underbrace{yy...x}}%
~\underset{j_{1}}{\underbrace{yx...x}}~y^{j_{2}}  \tag{2.5.}
\end{equation}%
and $y$ will commute with $x$ from $j_{1}$ times. Since we do this from $%
i_{2}$ times, we obtain the below formula for the function 
\begin{equation}
\alpha \left( i,j\right) :\bigskip \{%
\begin{array}{c}
\xi ^{i_{2}j_{1}},\text{ if }i_{1}+j_{1}<n\text{ and }i_{2}+j_{2}<n \\ 
\xi ^{i_{2}j_{1}}a,~\text{if }i_{1}+j_{1}\geq n\text{ and }i_{2}+j_{2}<n \\ 
\xi ^{i_{2}j_{1}}b,\text{if }i_{1}+j_{1}<n\text{ and }i_{2}+j_{2}\geq n \\ 
\xi ^{i_{2}j_{1}}ab,\text{if }i_{1}+j_{1}\geq n\text{ and }i_{2}+j_{2}\geq n%
\end{array}
\tag{2.6.}
\end{equation}

\textbf{The algorithm\medskip }

\textbf{Input:} $n,e_{i},$ $e_{j},i_{1},i_{2},j_{1},j_{2}$

\textit{Step 1.} Find $n-$ary decomposition $i_{n}$ and $j_{n}$ for the
numbers $i$ and $j.$

\textit{Step 2.} Compute $i_{n}\circ j_{n}$ in the group $\mathbb{Z}%
_{n}^{n}. $

\textit{Step 3.} Compute $\alpha \left( i,j\right) $ using formula $\left(
2.6\right) .$

\textbf{Output:} $e_{i}e_{j}$%
\begin{equation*}
\end{equation*}

\textbf{3. Octonion algebra in a symbol algebra of degree three} 
\begin{equation*}
\end{equation*}

In the following, we will \ show what is the connexion between the octonion
algebras, algebras of dimension $8,$ and symbol algebras of degree $3,$
proving that in all symbol algebra of degree three we can find an octonion
algebra without division.

Let $S$ be an associative algebra of degree three. For $z\in S,$ let $%
P\left( X,z\right) $ be the characteristic polynomial for the element $z$%
\begin{equation}
P\left( X,z\right) =X^{3}-\tau \left( z\right) X^{2}+\pi \left( z\right)
X-\eta \left( z\right) \cdot 1,  \tag{3.1.}
\end{equation}%
where $\tau $ is the linear form, $~\pi $ is the quadratic form and $\eta $
the cubic form.\medskip

\textbf{Proposition 3.1.} ([Fa; 88], Lemma) \textit{With the above
notations, denoting by} \textit{\ }$z^{\ast }=z^{2}-\tau \left( z\right)
z+\pi \left( z\right) \cdot 1,\,\,$\textit{for an associative algebra of
degree three,} \textit{we have:}

\textit{i) }$\pi \left( z\right) =\tau \left( z^{\ast }\right) .$

\textit{ii) }$z^{\ast \ast }=\eta \left( z\right) z.\Box \medskip \smallskip
\medskip $

An associative finite dimensional $K$-algebra $A$ is \textit{semisimple} if
it can be expressed as a finite and unique direct sum of simple algebras. An
associative $K$-algebra $A$ is \textit{separable} if for every field
extension $K\subset L$ the \ algebra $A\otimes _{K}L$ is semisimple. \ We
have that any central simple algebra is a separable algebra over its center
(see [Ha; 00], p.463). A \textit{Hurwitz} algebra $A$ is a unital (not
necessarily associative) algebra over $K$ together with a nondegenerate
quadratic form $n$ which satisfies $n\left( xy\right) =n\left( x\right)
n\left( y\right) ,x,y\in A.\medskip $

\textbf{Theorem 3.2.} ([Ja; 81], Theorem 6.2.3) \textit{Let }$A$ \textit{be
a finite-dimensional algebra with unity over the field} $K$ \textit{and} $%
\varphi :A\rightarrow K~\ $\textit{be a nondegenerate quadratic form such
that} $\varphi \left( xy\right) =\varphi \left( x\right) \varphi \left(
y\right) $ \textit{for all} $x,y\in A.$\textit{Then the algebra} $A$ \textit{%
has dimension} $1,2,4$ \textit{or} $8.$ \textit{If} $\dim A\in \{4.8\},$ $A$ 
\textit{is a quaternion \ or an octonion algebra}. $\Box \medskip \medskip $

Let $S=\left( \frac{a,b}{K,\xi }\right) $ be a symbol algebra of degree $n.$
For $n=3,$ the obtained symbol algebra has dimension $9$ over the field $K$
and, since an octonion algebra generalizes the quaternion algebra and has
dimension $8$ less than $9,$ we ask if we can find a relation between a
symbol algebra of degree three and an octonion algebra.\medskip

\textbf{Proposition 3.3.} ([Fa; 88], Theorem) \textit{If} $A$ \textit{is an
associative algebra of degree three over a field} $K$ \textit{containing the
cubic root of the unity, }$\xi ,$\textit{\ then, using notations from
Proposition 3.1, the quadratic form} $\pi $ \textit{permits compositions} $%
\pi \left( z\circ w\right) =\pi \left( z\right) \pi \left( w\right) $ 
\textit{on} $\widetilde{A}=\{u\in A~/~\tau \left( u\right) =0\}$ \textit{%
relative to the product} 
\begin{equation*}
z\circ w=\xi zw-\xi ^{2}wz-\frac{2\xi +1}{3}\tau \left( zw\right) \cdot
1,z,w\in \widetilde{A}.
\end{equation*}%
\textit{If} $A$ \textit{is separable over} $K$, \textit{therefore the
quadratic form }$\pi $\textit{\ is nondegenerate and} $\ $\textit{we can
find a new product "}$\bigtriangledown $\textit{" on} $(\widetilde{A},\circ
) $ \textit{such that }$\left( \widetilde{A},\bigtriangledown \right) $ 
\textit{is a Hurwitz algebra.}$\medskip \medskip \Box \medskip $

Since $S$ is separable over $K$ it results that $\pi $ is a nondegenerate
quadratic form on $S$ and it is also nondegenerate on $\widetilde{S}$, then
there is an element $u\in \widetilde{S}$ such that $\pi \left( u\right) \neq
0.$ Using some ideas given in [Ka; 53], let $v=\frac{u^{2}}{\pi \left(
u\right) }.\medskip $

\textbf{Proposition 3.4.} \textit{The linear maps }$R_{v}^{\circ }:%
\widetilde{S}\rightarrow \widetilde{S},R_{v}^{\circ }\left( x\right) =x\circ
v$ \textit{and} $L_{v}^{\circ }:\widetilde{S}\rightarrow \widetilde{S}%
,L_{v}^{\circ }\left( x\right) =x\circ v$ \textit{are bijective}$.\medskip
\medskip $

\textbf{Proof.} \ Let $R_{v}^{\circ }:\widetilde{S}\rightarrow \widetilde{S}%
,R_{v}^{\circ }\left( x\right) =x\circ v.$ Since $\pi \left( R_{v}^{\circ
}\left( x\right) \right) =\pi \left( x\circ v\right) =\pi \left( x\right)
\pi \left( v\right) =\pi \left( x\right) ,$ if $R_{v}^{\circ }\left(
x\right) =0$ it results that $\pi \left( x\right) =0,$ Using that $\pi $ is
nondegenerate, we obtain\textit{\ } $x=0,$ therefore \bigskip $R_{v}^{\circ
}~$is bijective. $\Box \medskip $

From the \ above proposition, on $\widetilde{S\text{ }},$ we define a new
multiplication 
\begin{equation*}
z\bigtriangledown w=(R_{v}^{\circ -1}\left( z\right) )\circ (L_{v}^{\circ
-1}\left( w\right) ),w,z\in \widetilde{S}.
\end{equation*}%
We have that $v\circ v$ is the unity element and $\pi \left(
z\bigtriangledown w\right) =\pi \left( z\right) \pi \left( w\right) $.
Indeed, it results\newline
$z\bigtriangledown \left( v\circ v\right) =(R_{v}^{\circ -1}\left( z\right)
)\circ (L_{v}^{\circ -1}\left( v\circ v\right) )=$\newline
$=$ $(R_{v}^{\circ -1}\left( z\right) )\circ L_{v}^{\circ -1}\left(
L_{v}^{\circ }\left( v\right) \right) =$\newline
$=(R_{v}^{\circ -1}\left( z\right) )\circ v=R_{v}^{\circ }\left(
R_{v}^{\circ -1}\left( z\right) \right) =z=$\newline
$=\left( v\circ v\right) \ \bigtriangledown z\ $and $\ $\newline
$\pi \left( z\bigtriangledown w\right) =\pi \left( (R_{v}^{\circ -1}\left(
z\right) )\circ (L_{v}^{\circ -1}\left( w\right) )\right) =$\newline
$=\pi \left( R_{v}^{\circ -1}\left( z\right) \right) \pi \left( L_{v}^{\circ
-1}\left( w\right) \right) =$\newline
$=\pi \left( z\right) \pi \left( w\right) ,$ since $\pi \left( z\right) =\pi
\left( R_{v}^{\circ }\left( R_{v}^{\circ -1}\left( z\right) \right) \right)
=\pi \left( R_{v}^{\circ -1}\left( z\right) \right) .$ Therefore the algebra 
$\left( \widetilde{S},\bigtriangledown \right) $ is a Hurwitz algebra of
dimension $8$ and, from Theorem 4.1, we obtain that $\left( \widetilde{S}%
,\bigtriangledown \right) $ is an octonion algebra with the norm $\pi $. An
octonion algebra $A$ with the norm $\pi $ is a division algebra if $\pi
\left( x\right) =0$ implies $x=0,$ for $x\in A.$ This algebra is a not a
division algebra since, from Proposition 3.1, we have \ $0=\tau \left( \eta
\left( x\right) x\right) =\tau \left( (x^{\ast })^{\ast }\right) =\pi \left(
x^{\ast }\right) ,$ for the element $x\in \widetilde{S}.$\medskip

From the above, we proved the following theorem:\medskip

\textbf{Theorem 3.5.} \textit{Let} $S=\left( \frac{a,b}{K,\xi }\right) $ 
\textit{be a symbol algebra of degree} $3.$ \textit{On the vector space} $%
\widetilde{S},$ \textit{we define the following products:}%
\begin{equation*}
z\circ w=\xi zw-\xi ^{2}wz-\frac{2\xi +1}{3}\tau \left( zw\right) \cdot
1,z,w\in \widetilde{S}
\end{equation*}%
\textit{and}%
\begin{equation*}
z\bigtriangledown w=(R_{v}^{\circ -1}\left( z\right) )\circ (L_{v}^{\circ
-1}\left( w\right) ),z,w\in \widetilde{S}.
\end{equation*}%
\textit{Therefore }$\left( \widetilde{S},\bigtriangledown \right) $ \textit{%
is an octonion non-division algebra. \medskip }$\Box \medskip \medskip
\medskip $

Since always we can find an octonion non-division algebra in a symbol
algebra of degree three, a natural question appears: \ if we can find some
conditions which can determine when a symbol algebra of degree three is with
division or not, or, more simple, if there are examples of division symbol
algebras of degree three. Such as conditions was given in [Fla; 12] for
symbol algebras of degrees $n=4k+2$ (but not for $n=3$ or $n=5),$ in which
the author found some trace form criteria to determine if a symbol algebra
is with division. In the mentioned paper, the author don't provide examples,
as will do in the next section.

\begin{equation*}
\end{equation*}

\textbf{4. Examples of division symbol algebras of degree }$3$ \textbf{and} $%
5$\textbf{\ }%
\begin{equation*}
\end{equation*}

In this section, we determine certain class of non-division symbol algebras
using Theorem 1.1 and some properties of ramification theory in algebraic
number fields, for example the decomposition of a prime ideal in $\mathbb{Z}$
in the ring of integers of a cyclotomic field or the decomposition of a
prime ideal in the ring $\mathbb{Z}\left[ \xi \right] $ in the ring of
integers of a Kummer field (see Theorem 1.2.and Theorem 1.3.). We also
provide examples of division symbol algebras of degree $3$ and $5.\medskip $

\smallskip

\textbf{Proposition 4.1.}\textit{Let} $\epsilon $ \textit{be a primitive
root of order} $3$ \textit{of unity and let }$K=\mathbb{Q}\left( \epsilon
\right) $\textit{\ be the cyclotomic field.} \textit{Let} $\alpha \in
K^{\ast },$ $p$ \textit{a prime rational integers,} $p\neq 3$ \textit{and
let the Kummer field} $L=K\left( \sqrt[3]{\alpha }\right) $ \textit{such that%
} $\alpha $ \textit{is a cubic residue modulo} $p.$ \textit{Let} $h_{L}$ 
\textit{be the class number of} $L.$ \textit{Then, the symbol algebras} $%
A=\left( \frac{\alpha ,p^{h_{L}}}{K,\epsilon }\right) $\textit{\ is
non-division.}

\smallskip 

\textbf{Proof.} Since $p$ is a prime rational integer, $p\neq 3,$ it results 
$p\equiv 1$ (mod $3$) or $p\equiv 2$ (mod $3$).\newline
\textbf{Case 1:} is $p\equiv 2$ (mod $3$).\newline
We know that the ring of integers of $K$ is $\mathbb{Z}\left[ \epsilon %
\right] $ and it is a principal ring. According to Theorem 1.2 it results
that $p$ is inert in the ring $\mathbb{Z}\left[ \epsilon \right] .$ If we
denote with $\mathcal{O}_{L}$ the ring of integers of the Kummer field $L$
and knowing that the cubic residual symbol $\left( \frac{\alpha }{p}\right)
_{3}=1,$ we apply Theorem 1.3 and we obtain that: 
\begin{equation*}
p\mathcal{O}_{L}=P_{1}P_{2}P_{3},
\end{equation*}%
where $P_{1},$ $P_{2},$ $P_{3}$ are conjugate prime ideals. We obtain that: 
\begin{equation*}
\left( p\mathcal{O}_{L}\right)
^{h_{L}}=P_{1}^{h_{L}}P_{2}^{h_{L}}P_{3}^{h_{L}}.
\end{equation*}%
Therefore, there exists a principal ideal $I$ in the ring $\mathcal{O}_{L}$
such that $\left( p\mathcal{O}_{L}\right) ^{h_{L}}=N_{L/K}\left( I\right) .$
It results that there exists $x$$\in \mathcal{O}_{L}$ such that $%
p^{h_{L}}=N_{L/K}\left( x\right) .$ Applying Theorem 1.1 and Remark 1.6, we
obtain that the symbol algebras $A=\left( \frac{\alpha ,p^{h_{L}}}{%
K,\epsilon }\right) $ is non-division.\newline
\textbf{Case 2:} is $p\equiv 1$ (mod $3$).\newline
Applying Theorem 1.2 it results that 
\begin{equation*}
p\mathbb{Z}\left[ \epsilon \right] =p_{1}\mathbb{Z}\left[ \epsilon \right]
p_{2}\mathbb{Z}\left[ \epsilon \right] ,
\end{equation*}%
where $p_{1},p_{2}$$\in $$\mathbb{Z}\left[ \epsilon \right] .$\newline
Since $\alpha $ is a cubic residue modulo $p$ we obtain that the cubic
residual symbols $\left( \frac{\alpha }{p_{1}}\right) _{3}=\left( \frac{%
\alpha }{p_{2}}\right) _{3}=1.$ Applying Theorem 1.3 it results that 
\begin{equation*}
p\mathcal{O}_{L}=p_{1}\mathcal{O}_{L}p_{2}\mathcal{O}%
_{L}=P_{11}P_{12}P_{21}P_{22}P_{31}P_{32},
\end{equation*}%
where $P_{i1}$ and $P_{i2},$ $i=\overline{1,3}$ are conjugate prime ideals.
We obtain that: 
\begin{equation*}
\left( p\mathcal{O}_{L}\right) ^{h_{L}}=\left( P_{11}P_{21}\right)
^{h_{L}}\left( P_{12}P_{22}\right) ^{h_{L}}\left( P_{13}P_{23}\right)
^{h_{L}}.
\end{equation*}%
From this, as in Case 1, we obtain that the symbol algebras $A=\left( \frac{%
\alpha ,p^{h_{L}}}{K,\epsilon }\right) $ is non-division.\textit{\medskip }$%
\Box $
\smallskip

\textbf{Corollary 4.2.} \textit{Let} $q$ \textit{be an odd prime positive
integer and} $\xi $ \textit{be a primitive root of order} $q$ \textit{of
unity and let }$K=\mathbb{Q}\left( \xi \right) $\textit{\ be the cyclotomic
field.} \textit{Let} $\alpha \in K^{\ast },$ $p$ \textit{a prime rational
integers,} $p\neq 3$ \textit{and let the Kummer field} $L=K\left( \sqrt[q]{%
\alpha }\right) $ \textit{such that} $\alpha $ \textit{is a} $q$ \textit{%
power residue modulo} $p.$ \textit{Let} $h_{L}$ \textit{be the class number
of} $L.$ \textit{Then, the symbol algebras} $A=\left( \frac{\alpha ,p^{h_{L}}%
}{K,\xi }\right) $\textit{\ is non-division.}\newline
\smallskip \newline
\textbf{Proof.} The proof is similar with the proof of Proposition 4.1.%
\textit{\medskip }$\Box $%
\begin{equation*}
\end{equation*}%
\smallskip \newline
In the following, we will give some examples of division symbol algebras.
Using the computer algebra system MAGMA, we found some examples of division
symbol algebras of degree $3$ and of degree $5.$

\medskip

\textbf{Example 4.3.} \smallskip 
Q:=Rationals();\newline
E:=CyclotomicField($3$);\newline
a:=RootOfUnity($3$);\newline
a;\newline
$Et<t>:=$PolynomialRing(E);\newline
E;\newline
f := $t^{3}-7;$\newline
$K<b>:=$NumberField(f);\newline
K;\newline
$b^{3};$\newline
NormEquation(K,$11$);\newline
NormEquation(K,$11$+a);\newline
NormEquation(K,$11^{3}$);\newline
NormEquation(K,$(11+a)^{3}$);\newline
NormEquation(K,$5$);\newline
NormEquation(K,$5$+a);\newline
NormEquation(K,$5^{3}$);\newline
NormEquation(K,$(5+a)^{3}$);\newline
Evaluate\newline
$zeta_{3}$\newline
Cyclotomic Field of order $3$ and degree $2$\newline
Number Field with defining polynomial $t^{3}-7$ over E\newline
$7$\newline
false\newline
false\newline
true [$-11\ast zeta_{3}-11$]\newline
true [$10\ast zeta_{3}-1$]\newline
false\newline
false\newline
true [$5\ast zeta_{3}$]\newline
true [$zeta_{3}+5$]\newline
\begin{equation*}
\end{equation*}%
\smallskip \newline
\textbf{Example 4.4.}\newline
\newline
Q:=Rationals();\newline
F:=CyclotomicField($5$);\newline
a:=RootOfUnity($5$);\newline
a;\newline
$Ft<t>:=$PolynomialRing(F);\newline
F;\newline
f := $t^{5}-13;$\newline
$K<b>:=$NumberField(f);\newline
K;\newline
$b^{3};$\newline
NormEquation(K,$11$);\newline
NormEquation(K,$11$+a);\newline
Evaluate\newline
$zeta_{5}$\newline
Cyclotomic Field of order $5$ and degree $4$\newline
Number Field with defining polynomial $t^{3}-13$ over F\newline
$13$\newline
false\newline
false\newline
\medskip
From the above examples, using Teorem 1.1 and Remark 1.6, we obtain that the
symbol algebras $\left( \frac{7,11}{\mathbb{Q}\left( \epsilon \right)
,\epsilon }\right) ,$ $\left( \frac{7,11+\epsilon }{\mathbb{Q}\left(
\epsilon \right) ,\epsilon }\right) ,$ $\left( \frac{7,5}{\mathbb{Q}\left(
\epsilon \right) ,\epsilon }\right) ,$ $\left( \frac{7,5+\epsilon }{\mathbb{Q%
}\left( \epsilon \right) ,\epsilon }\right) $ are division symbol algebras
of degree $3$ and $\left( \frac{13,11}{\mathbb{Q}\left( \xi \right) ,\xi }%
\right) ,$ $\left( \frac{13,11+\xi }{\mathbb{Q}\left( \xi \right) ,\xi }%
\right) $ are division symbol algebras of degree $5.$

\smallskip 

In the following we determine some split symbol algebras.\newline
\smallskip \newline
\textbf{Proposition 4.5.} \textit{Let} $n\geq 2$ \textit{be an arbitrary
positive integer and let} $\xi $ \textit{be a primitive root of order} $n$ 
\textit{of unity.} \textit{Let} $K$ \textit{be a finite field whose} $%
char(K) $ \textit{does not divide} $n,$ $a,b$ $\in K^{\ast }$ \textit{and let%
} $S$ \textit{be the symbol algebra} $S=\left( \frac{a,~b}{K\left( \xi
\right) ,\xi }\right) .$ \textit{Then} $S$ \textit{is a split algebra.}%
\newline
\smallskip \newline
\textbf{Proof.} Applying Theorem 1.5 and Theorem 1.1 we obtain that $S$ is a
split algebra.\newline
Another solution is to apply Theorem 1.4 and Remark 1.6. \textit{\medskip }$%
\Box \medskip $ \newline
\textbf{Remark 4.6.} For $n$ a prime number, all symbol algebras from the
above proposition are non-division algebras.%
\[\]
\smallskip \newline
\begin{equation*}
\end{equation*}

\textbf{Conclusions.} In this paper, we gave an algorithm for compute
quickly the elements from the basis in a symbol algebra of degree $n$ and we
find an octonion non-division algebra in a symbol algebra of degree three$.$
We also provide some examples of division symbol algebras of degree three \
and five. Starting from results obtained in this paper, we intend to find,
in a furter research, more conditions for a symbol algebra of degree $n$ to
be with division. 
\begin{equation*}
\end{equation*}%
\textbf{Acknowledgements.} The second author thanks Professor Tamas Szamuely
for helpful discussions on this topic. This paper \ is supported by the
grants of CNCS (Romanian National Council of Research):
PN-II-ID-WE-2012-4-169,\ for the first author and PN-II-ID-WE - 2012 - 4
-161, for the second author. 
\begin{equation*}
\end{equation*}

\textbf{References}%
\begin{equation*}
\end{equation*}

[Al, Io; 84] Albu, T. Ion I.D., \textit{Chapters of the algebraic Number
Theory} (in Romanian), Ed. Academiei, Bucharest, 1984.

[Ba; 09] Bales, J. W., \ \textit{A Tree for Computing the Cayley-Dickson
Twist}, \ Missouri J. Math. Sci., \textbf{21(2)}(2009), 83-93.

[Fa; 88]$~$\ J.R. Faulkner, \textit{Finding octonion algebras in associative
algebras,} Proc. Amer. Math. Soc. \textbf{104}(4)(1988), 1027-1030.

[Fla; 12] R. Flatley, \textit{Trace forms of Symbol Algebras}, Algebra
Colloquium, \textbf{19}(2012), 1117-1124.

[Gi, Sz; 06] Gille, P., Szamuely, T., \textit{Central Simple Algebras and
Galois Cohomology}, Cambridge University Press, 2006.

[Ha; 00] M. Hazewinkel, \textit{Handbook of Algebra}, Vol. 2, North Holland,
Amsterdam, 2000.

[Ir, Ro; 92] Ireland,K., Rosen M. \textit{A Classical Introduction to Modern
Number Theory}, Springer Verlag, 1992.

[Jac; 81] N. Jacobson, \textit{Structure theory of Jordan algebras}, Lecture
Notes in Mathematics, Vol. 5, The University of Arkansas, Fayetteville, 1981.

[Ja, Ya; 13]M. Jafari, Y. Yayli,\textit{Rotation in four dimensions via
Generalized Hamilton operators Kuwait Journal of Science,} vol 40 (1) June
2013, p.67-79.

[Jan; 73]G.J. Janusz, \textit{Algebraic number fields}, Academic Press,
London, 1973.

[Fl, Sa; 13] C. Flaut, D. Savin, \textit{Some properties of the symbol
algebras of degree} $3$, submitted.

[Ka; 53] L. Kaplansky, \textit{Infinite-dimensional quadratic forms
admitting composition, }Proc. Amer. Math., \textbf{\ 4(}1953), 956-960.

[Led; 05] A. Ledet, \textit{Brauer Type Embedding Problems }, American
Mathematical Society, 2005.

[Lem; 00] F. Lemmermeyer, \textit{Reciprocity laws, from Euler to Eisenstein 
}, Springer-Verlag, Heidelberg, 2000.

[Mil; 10] J.S. Milne, \textit{Class Field Theory},
http://www.math.lsa.umich.edu/~jmilne.

[Mi; 71] \ J. Milnor, \ \textit{Introduction to Algebraic K-Theory, }Annals
of Mathematics Studies, Princeton Univ. Press, 1971.

[Pi; 82] Pierce, R.S., \textit{Associative Algebras}, Springer Verlag, 1982.

[Sa, Fl, Ci; 09] D. Savin, C. Flaut, C. Ciobanu, \textit{Some properties of
the symbol algebras}, Carpathian Journal of Mathematics , vol. 25, No. 2
(2009), p. 239-245.

[Sc; 66] Schafer, R. D., \textit{An Introduction to Nonassociative Algebras,}
Academic Press, New-York, 1966.

[Sc; 54] Schafer, R. D., \textit{On the algebras formed by the
Cayley-Dickson process,} Amer. J. Math., \textbf{76}(1954), 435-446.%
\begin{equation*}
\end{equation*}

Cristina FLAUT

{\small Faculty of Mathematics and Computer Science,}

{\small Ovidius University,}

{\small Bd. Mamaia 124, 900527, CONSTANTA,}

{\small ROMANIA}

{\small http://cristinaflaut.wikispaces.com/}

{\small http://www.univ-ovidius.ro/math/}

{\small e-mail:}

{\small cflaut@univ-ovidius.ro}

{\small cristina\_flaut@yahoo.com}

\begin{equation*}
\end{equation*}

Diana SAVIN

{\small Faculty of Mathematics and Computer Science,}

{\small Ovidius University,}

{\small Bd. Mamaia 124, 900527, CONSTANTA,}

{\small ROMANIA}

{\small e-mail:} {\small savin.diana@univ-ovidius.ro}

{\small dianet72@yahoo.com}

\end{document}